\input amstex.tex
\input amsppt.sty

\define\E{\Bbb E}
\define\N{\Bbb N}
\define\R{\Bbb R}
\redefine\H{\Bbb H}
\define\C{\Bbb C}
\define\Z{\Bbb Z}
\define\tr{\operatorname{Tr}}

\NoBlackBoxes
\NoRunningHeads

\topmatter
\title CLT for spectra of submatrices of Wigner random matrices II. 
Stochastic evolution
\endtitle

\author Alexei Borodin
\endauthor

\abstract We show that the global fluctuations of spectra of GOE and GUE
matrices and their principal submatrices executing Dyson's Brownian motion
are Gaussian in the limit of large matrix dimensions. For nested submatrices
one obtains a limiting three-dimensional generalized Gaussian process; its
restrictions to two-dimensional sections that are monotone in matrix sizes
and time moments coincide with the two-dimensional Gaussian Free Field with 
zero boundary conditions.
The proof is by moment convergence, and it extends to more general 
Wigner matrices and their stochastic evolution.   
\endabstract
\endtopmatter

\document

\subhead Introduction
\endsubhead The fact that the global spectral fluctuations of a GOE or a GUE
random matrix evolving under Dyson's Brownian Motion, are asymptotically
Gaussian is well-known, see \S4.3.3 in \cite{AGZ} and references therein, and also
\cite{S} for a general $\beta$ analog. On the other hand, it was shown in \cite{B}
that the global fluctuations of spectra of various principal submatrices of a single 
GOE or GUE matrix are also Gaussian. The goal of this note is to 
put these two statements together. 

We prove the asymptotic Gaussian behavior for submatrices of a class of stochastically
evolving Wigner random matrices that includes Dyson's Brownian Motion for GOE and GUE. 
The proof is by the method of moments, and the argument is slightly more general
than the one presented in \cite{AGZ} for a single Wigner matrix. 

We also compute the resulting covariance kernel explicitly. In the case of nesting
submatrices, it represents a three-dimensional generalized Gaussian process, where 
one dimension comes from the position of the spectral variable, the second dimension 
reflects the size of the submatrix, and the third dimension is the time variable. 
When restricted to the two-dimensional sections that are monotone in matrix size and 
time variables, it reproduces the two-dimensional Gaussian Free Field (GFF) with zero
boundary conditions.

In the case of GUE, the appearance of GFF on monotone sections 
could have been predicted from the determinantal sctructure of the correlation
functions \cite{FF}, \cite{ANvM1}, and from the analysis of \cite{BF} that showed how
such a structure leads to GFF covariances in the global asymptotic regime. However,
the complete three-dimensional covariance structure seems to be unaccessible via
that approach for example because the spectra of the full set of submatrices evolve in 
a non-Markovian
way \cite{ANvM2}.

\noindent{\it Acknowledgements.\/} This work was partially supported by NSF grant
DMS-1056390.

\subhead Wigner matrices
\endsubhead
Let $\{Z_{ij}(t)\}_{j> i\ge 1,t\in\R}$ and $\{Y_i(t)\}_{i\ge 1,t\in\R}$ be two families 
of independent identically distributed real-valued stochastic (not necessarily Markov) 
processes with zero mean 
such that for any $k\ge 1$
$$
\max_{t\in\R}(\,\E |Z_{12}(t)|^k,\E|Y_1(t)|^k)<\infty.
$$

Set $c(s,t)=\frac 12\E Y_1(s) Y_1(t)$ and assume that
$$
\gathered
c(s,t)\ge 0, \qquad c(t,t)\equiv 1,\\ 
\E Z_{12}(s)Z_{12}(t)\equiv c(s,t),\qquad
\E Z_{12}^2(s)Z_{12}^2(t)\equiv 2c(s,t)^2+1.
\endgathered
$$
Note that by Cauchy's inequality $c(s,t)\le \sqrt{c(s,s)c(t,t)}=1$. We say that a
function $c(s,t)$ is admissible if it arises in this way.  

One possibility for the above relations to be satisfied is to take all 
$\{2^{-\frac 12} Y_i(t)\}$ and $\{Z_{ij}(t)\}$ to be independent standard 
Ornstein-Uhlenbeck processes on $\R$; then $c(s,t)=\exp(-|s-t|)$. We will refer
to this possibility as to {\it Gaussian specialization\/}. 

Define a (real symmetric) {\it time-dependent Wigner matrix\/} $X(t)=[X(i,j\mid
t)]_{i,j\ge 1}$ by
$$
X(i,j\mid t)=X(j,i\mid t)=\cases Z_{ij}(t),& i<j,\\ Y_i,&i=j.\endcases
$$

An Hermitian variation of the same definiton is as follows:
Let $\{Z_{ij}\}_{j> i\ge 1}$ now be complex-valued (i.i.d. mean zero) stochastic
processes with the same uniform bound on all
moments. Denote $d(s,t)=\E Y_1(s) Y_1(t)$ and assume that
$$
\gathered
d(s,t)\ge 0, \qquad d(t,t)\equiv 1,\\ \E Z_{12}(s)Z_{12}(t)\equiv 0,\quad 
\E Z_{12}(s)\overline{Z_{12}(t)}\equiv d(s,t),\quad \E
|Z_{12}(s)|^2|Z_{12}(t)|^2\equiv d(s,t)^2+1.
\endgathered
$$
We will also say that a function $d(s,t)$ is admissible if it arises in this way.  

There is also a Gaussian specialization that corresponds to $\{Y_i(t)\}$ and 
$\{2^{\frac 12}\Re Z_{ij}(t)\}$, $\{2^{\frac 12}\Im Z_{ij}(t)\}$ being independent
standard Ornstein-Uhlenbeck processes on $\R$; in that case $d(s,t)=\exp(-|s-t|)$.

Define an {\it Hermitian\/} time-dependent Wigner matrix $X(t)=[X(i,j\mid
t)]_{i,j\ge 1}$ by
$$
X(i,j\mid t)=\overline{X(j,i\mid t)}=\cases Z_{ij}(t),& i<j,\\ Y_i,&i=j.\endcases
$$

Under the Gaussian specializations, the matrix stochastic processes defined above
are called {\it Dyson's Brownian motions\/}. Traditionally one distinguishes the 
two cases by a parameter $\beta$ that takes value 1 in the the real symmetric case
and value 2 in the Hermitian case. The random matrices arising at a single time
moment are said to belong to the
Gaussian Orthogonal Ensemble (GOE) in the $\beta=1$ case, and Gaussian Unitary
Ensemble (GUE) in the $\beta=2$ case.

\subhead The height function
\endsubhead
For any finite set $B\subset \{1,2,\dots\}$ denote by $X_B$ the $|B|\times |B|$
submatrix of a matrix $X$ formed by the intersections of the rows and columns of 
$X$ marked by elements of $B$. 

The {\it height function\/} $H$ associated to a time-dependent Wigner matrix
$X$ is a random integer-valued function on $\R\times \R_{\ge 1}\times\R$ defined by
$$
H(x,y,t)=\sqrt{\frac{\beta\pi}{2}}\,\bigl\{\text{the number of eigenvalues of
}X_{\{1,2,\dots,{[y]}\}}(t)\text{ that
are }\ge x\bigr\}.
$$

More generally, let $A=\{a_n\}_{n\ge 1}$ be an arbitrary sequence of pairwise 
distinct natural numbers. Then we define the height function $H_A$ via 
$$
H_A(x,y)=\sqrt{\frac{\beta\pi}{2}}\,\bigl\{\text{the number of eigenvalues of
}X_{\{a_1,\dots,a_{[y]}\}}(t)\text{ that
are }\ge x\bigr\}.
$$
The first definition corresponds to $A=\N$. 

The convenience of the constant prefactor $\sqrt{{\beta\pi}/{2}}$ will be
evident shortly.

\subhead A three-dimensional Gaussian field
\endsubhead
Let $c(s,t)$ be an admissible function as defined above. 
Set $\H=\{z\in\C\mid \Im z>0\}$ and
introduce a function $$
C:(\H\times\R)\times(\H\times\R)\to \R\cup\{+\infty\}
$$ 
via
$$
\multline
C(z,s;w,t)=\frac 1{2\pi}\ln\left|\frac{c(s,t)\min(|z|^2,|w|^2)-zw}
{c(s,t)\min(|z|^2,|w|^2)-z\overline{w}}\right|\\ =
\cases  -\frac 1{2\pi}\ln\left|\dfrac{c(s,t) z-w}
{c(s,t) z-\overline{w}}\right|,&|z|\le |w|,\\
-\frac 1{2\pi}\ln\left|\dfrac{c(s,t)w-z}
{c(s,t)w-\overline{z}}\right|,&|z|>|w|.
\endcases
\endmultline
$$
It is easy to see that for any $(s,t)$ with $c(s,t)<1$, $C(\,\cdot\,,s;\,\cdot\,,t)$ 
is a continuous function on $\H\times\H$. 
Note also that if $c(s,t)=1$ then
$$
C(z,s;w,t)=-\frac 1{2\pi}\ln\left|\frac{z-w}
{z-\overline{w}}\right|
$$
is the Green function for the Laplace operator on $\H$ with Dirichlet boundary
conditions. Viewed as a function in $(z,w)$, it represents the covariance for the
two-dimensional Gaussian Free Field on $\H$ with zero boundary conditions.

\proclaim{Proposition 1} For any admissible function $c(s,t)$ as above, there exists a
generalized Gaussian process on $\H\times\R$ with the covaraince kernel $C(z,s;w,t)$ as 
above. More exactly, for any finite family of test functions $f_m(z)\in C_0(\H\times \R)$ the covariance matrix
$$
cov(f_k,f_l)=\int_\H\int_\H f_k(z,s)f_l(w,t) C(z,s;w,t) \,dzd\bar{z}ds\,dwd\bar{w}dt,\quad
k,l=1,\dots,M,
$$
is positive-definite.
\endproclaim

Denote the resulting generalized Gaussian process by $\Cal G_{c(s,t)}$.

A proof of Propositon 1 will be given later.

\subhead Complex structure
\endsubhead
Let $A$ be a sequence of pairwise distinct integers. The height function $H_A(x,y,t)$ 
(or $H(x,y,t)=H_\N(x,y,t)$)
is naturally defined on $\R\times \R_{\ge 1}\times\R$. Having the large parameter $L$, we would
like to scale $(x,y)\mapsto (L^{-\frac 12}x,L^{-1}y)$, which lands us in $\R\times
\R_{> 0}\times\R$.

Wigner's semicircle law implies that for any $t\in\R$, 
with $L\gg 1$, $x\sim L^\frac 12$, $y\sim L$,
after rescaling with overwhelming probability
the eigenvalues (or, equivalently, the places of growth of the height function in $x$-direction) 
are
concentrated in the domain
$$
\bigl\{(x,y)\in\R\times\R_{>0}\mid -2\sqrt{y}\le x\le 2\sqrt{y}\bigr\}.
$$

Let us identify the interior of this domain with $\H$ via the map
$$
\Omega:(x,y)\mapsto\frac x2+ i\sqrt{y-\left(\frac x2\right)^2}.
$$
Its inverse has the form
$$
\Omega^{-1}(z)=(x(z),y(z))=(2\Re(z),|z|^2).
$$
Note that this map sends the boundary of the domain to the real line.

Thanks to $\Omega$ we can now speak of the height function $H_A$ as being defined on
$\H\times \R$; we will use the notation
$$
H_A^\Omega(z;t)=H_A(L^{\frac 12}x(z), Ly(z),t),\qquad z\in\H.
$$
Note that we have incorporated rescaling in this definition.

\subhead{Main result}\endsubhead Let $X$ be a (real symmetric or Hermitian)
time-dependent Wigner matrix.
We argue that the centralized random height function
$$
H^\Omega(z;t)-\E H^\Omega(z;t),\qquad z\in\H,\quad t\in\R,
$$
viewed as distribution, converges as $L\to\infty$ to the generalized Gaussian process
$\Cal G_{c(s,t)}$ with $c(s,t)=\frac{\beta}2\,\E Y_1(s)Y_1(t)$.

One needs to verify the convergence on a suitable set of test functions.
The exact statement that we prove is the following.

\proclaim{Theorem 2}
Pick $\tau\in \R$, $y>0$, and $k\in \Z_{\ge 0}$. Define a moment of the random
height function by
$$
M_{\tau,y,k}=\int_{-\infty}^{+\infty} x^k \bigl(H(L^{\frac 12}x,Ly,\tau)
-\E H(L^{\frac
12}x,Ly,\tau)\bigr)dx.
$$
Then as $L\to\infty$, these moments converge, in the sense of finite dimensional
distributions, to the moments of $\Cal G_{c(s,t)}$ defined as
$$
\Cal M_{\tau,y,k} = \int_{z\in\H,|z|^2=y} (x(z))^k\, \Cal G_{c(s,t)}(z;\tau)\,
\frac{dx(z)}{dz}\,dz.
$$
\endproclaim

\subhead{Monotone sections as two-dimensional Gaussian Free Fields}\endsubhead Consider
a time-dependent Wigner matrix and assume that the function $c(s,t)=\frac{\beta}2\,\E
Y_1(s)Y_1(t)$ is continuous and that it
has the following monotonicity property: For any $s\in\R$, $c(s,t)$ is
strictly increasing in $t\in (-\infty,s]$ and it is strictly decreasing in 
$t\in[s,+\infty)$. 
In other words, as time distance between matrices grows, the correlation decays. 
Further, assume that $c(s,t)\ne 0$ for any $s,t\in\R$. 

Let $\phi:\R\to\R_{>0}$ and $\psi:\R\to\R$ be a continuous nonincreasing and a 
continuous nondecreasing functions, and assume that for at least one of 
these functions the monotonicity is strict. 

Our goal is to consider the joint fluctuations of spectra of matrices 
$$
X_{\{1,2,\dots,[L\phi(t)]\}}(\psi(t)),\qquad t\in\R,
\tag 1
$$  
where $L\gg 1$ is a large parameter. By Wigner's semicircle law, the spectrum of such a
matrix scaled by $L^\frac12$ is concentrated on $[-2\sqrt{\phi(t)},2\sqrt{\phi(t)}]$.

The two extreme cases are
$\phi(t)\equiv\operatorname{const}$ (the size of the matrices is fixed and the
time is moving) and $\psi(t)\equiv \operatorname{const}$ (the time moment is fixed and
the size of the matrices is changing).

Let us choose a reference time moment $t_0\in\R$ and introduce a map
$$
\Xi:\{(x,t)\in \R\times \R\mid -2\sqrt{\phi(t)}<x <
2\sqrt{\phi(t)}\}\to\H
$$
as
$$
\Xi(x,t)=\cases c(\psi(t_0),\psi(t))\left(\dfrac x2+ i\sqrt{\phi(t)-\left(\dfrac
x2\right)^2}\right),&t\ge t_0,\\ \dfrac 1{c(\psi(t),\psi(t_0))}
\left(\dfrac x2+ i\sqrt{\phi(t)-\left(\dfrac x2\right)^2}\right), & t<t_0.
\endcases
$$
The continuity and monotonicity assumptions on $c,\phi,$ and $\psi$ 
are needed for $\Xi$ to be a bijection. Hence, its 
inverse is correctly defined, denote it as  $\Xi^{-1}(\zeta)=(x(\zeta),t(\zeta))$.

We can now view the height function $H$ for matrices \thetag{1} 
as a function on $\H$ via
$$
H^\Xi(\zeta)=H\bigl(L^\frac 12\cdot x(\zeta), L\cdot \phi(t(\zeta)), 
\psi(t(\zeta))\bigr). 
$$

Our main result implies that the centralized height function
$$
H^\Xi(\zeta)-\E H^\Xi(\zeta),\qquad \zeta\in\H,
$$ 
viewed as a distribution, converges as $L\to\infty$ to the Gaussian Free Field on 
$\H$ in the sense of Theorem 2.

\subhead Moments as traces
\endsubhead
Let us rescale the variable $x=L^{-\frac 12}u$ in the definition of $M_{\tau,y,k}$
and then integrate by parts. Since the derivative of the height function
$H(u,[Ly],t)$ in $u$ is
$$
\frac d{du} H(u,[Ly],t)=-\sqrt{\frac{\beta\pi}2} \sum_{s=1}^{[Ly]}
\delta(u-\lambda_s),
$$
where $\{\lambda_s\}_{1\le s\le[Ly]}$ are the eigenvalues of $X_{\{1,\dots,[Ly]\}}(t)$, 
we obtain
$$
\multline
M_{\tau,y,k}={L^{-\frac{k+1}2}} \sqrt{\frac{\beta\pi}2} \left(\sum_{s=1}^{[Ly]}
\frac{\lambda_s^{k+1}}{k+1}-\E\sum_{s=1}^{[Ly]}
\frac{\lambda_s^{k+1}}{k+1}\right)\\=
\frac{L^{-\frac{k+1}2}} {k+1}\sqrt{\frac{\beta\pi}2}
\left(\tr\bigl(X_{\{1,\dots,[Ly]\}}^{k+1}(t)\bigr)-
\E \tr\bigl(X_{\{1,\dots,[Ly]\}}^{k+1}(t)\bigr)\right).
\endmultline
$$

We can now reformulate the statement of Theorem 2 as follows.

\proclaim{Theorem 2'} Let $X(t)$ be a time-dependent (real-symmetric or Hermitian)
Wigner matrix with
$c(s,t)=\frac \beta2 \E Y_1(s)Y_1(t)$. Let $k_1,\dots,k_m\ge 1$ be integers
and $y_1,\dots,y_m\in\R_{>0}$, $t_1,\dots,t_m\in\R$. 
Then the $m$-dimensional random vector
$$
\left(L^{-\frac {k_p}2}\biggl(\tr\bigl(X_{\{1,\dots,[Ly_p]\}}^{k_p}(t_p)\bigr)-
\E \tr\bigl(X_{\{1,\dots,[Ly_p]\}}^{k_p}(t_p)\bigr)\biggr)\right)_{p=1}^m
$$
converges (in distribution and with all moments) to the zero mean $m$-dimensional 
Gaussian random vector $(\xi_p)_{p=1}^m$ with covariance
$$
\multline
\E \xi_p\xi_q=\frac{2k_pk_q}{\beta \pi}
\oint\limits_{\Sb |z|^2=b_p\\ \Im z>0\endSb} \oint\limits_{\Sb |w|^2=b_q\\ \Im
w>0\endSb}
(x(z))^{k_p-1}(x(w))^{k_q-1}\\ \times
\frac 1{2\pi}\ln\left|\frac{c(t_p,t_q)\min(y_p,y_q)-zw}
{c(t_p,t_q)\min(y_p,y_q)-z\overline{w}}\right|\,
\frac{dx(z)}{dz}\frac{dx(w)}{dw}\,dzdw.
\endmultline
$$
\endproclaim

\subhead More general submatrices
\endsubhead In the spirit of \cite{B}, we will actually prove a more general claim
that involves arbitrary sequences of symmetric submatrices of the Wigner matrix 
that are sufficiently well-behaved. The exact statement is as follows. 

\proclaim{Theorem 2''} Let $X(t)$ be a time-dependent (real-symmetric or Hermitian)
 Wigner matrix with
$c(s,t)=\frac \beta2 \E Y_1(s)Y_1(t)$. 
Let $k_1,\dots,k_m\ge 1$ be integers, $t_1,\dots,t_m\in\R$, 
and let $B_1,\dots,B_m$ be subsets of $\N$ dependent on the large parameter $L$
so that there exist limits
$$
b_p=\lim_{L\to\infty} \frac{|B_p|}{L}>0,\qquad
b_{pq}=\lim_{L\to\infty}\frac{|B_p\cap B_q|}{L}\,,\qquad p,q=1,\dots,m.
$$
Then the $m$-dimensional random vector
$$
\left(L^{-\frac {k_p}2}\biggl(\tr\bigl(X_{B_p}^{k_p}(t_p)\bigr)-
\E \tr\bigl(X_{B_p}^{k_p}(t_p)\bigr)\biggr)\right)_{p=1}^m
\tag 2
$$
converges (in distribution and with all moments) to the zero mean $m$-dimensional Gaussian
random variable $(\xi_p)_{p=1}^m$ with the covariance
$$
\multline
\E \xi_p\xi_q=\frac{2k_pk_q}{\beta \pi}
\oint\limits_{\Sb |z|^2=b_p\\ \Im z>0\endSb} \oint\limits_{\Sb |w|^2=b_q\\ \Im
w>0\endSb}
(x(z))^{k_p-1}(x(w))^{k_q-1}\\ \times
\frac 1{2\pi}\ln\left|\frac{c(t_p,t_q)b_{pq}-zw}
{c(t_p,t_q)b_{pq}-z\overline{w}}\right|\,
\frac{dx(z)}{dz}\frac{dx(w)}{dw}\,dzdw.
\endmultline
\tag 3
$$
\endproclaim

Theorem 2" can also be viewed as the moment convergence
of the centralized height function $H_A(x,y,t)$ to a limiting 
generalized Gaussian process but we do not give further details here. 
The static variant of this convergence is discussed in \cite{B}. 

\subhead Proof of Theorem 2''
\endsubhead
The argument closely follows that given in Section 2.1.7 of
\cite{AGZ} in the case of one set $B_j\equiv B$, and the proof of Theorem 2' in
\cite{B} in the static case. One proves the
convergence of moments, which is sufficient to also claim the convergence
in distribution for Gaussian limits.

Any joint moment of the coordinates of \thetag{2} is written as a finite combination of
contributions corresponding to suitably defined graphs that are in their turn associated
to words. This reduction is explained in in Section 2.1.7 of
\cite{AGZ}. The key fact in the real-symmetric case is that averages of products of 
powers of matrix elements that involve at least one matrix element with exponent 1
vanish. The time-dependent analog of this fact is that averages of products of powers
of matrix elements taken at different time moments that involve one matrix element with
exponent 1 {\it at only one time moment\/} vanish. This clearly holds by independence of
matrix elements and our zero mean assumption. In the static Hermitian case, one needs in addition that $\E Z_{12}^2=0$. The time-dependent 
analog reads $\E Z_{12}(s)Z_{12}(t)=0$ for any $s,t\in\R$, which is one of our
assumptions. This allows the exact same reduction to go through in the time-dependent
setting. 

The only difference of the multi-set case from the one-set case is that one
needs to keep track of the {\it alphabets\/} the words are built from: A word
corresponding to coordinate number $p$ of \thetag{2} would have to be built from the
alphabet that coincides with the set $B_p$. Equivalently, the corresponding graphs
will have their vertices labeled by elements of $B_p$.

Since all sizes $|B_p|$ have order $L$, 
and $|B_1\cup\dots\cup B_m|=O(L)$, and also the moments of matrix elements at 
all times are uniformly bounded, the estimate
showing that all contributions not coming from matchings are negligible (Lemma 2.1.34
in \cite{AGZ}) carries over without difficulty. It only remains to compute the
covariance.

For real symmetric Wigner matrices in the one-set case the limits of the
variances of the coordinates of \thetag{2} are given by (2.1.44) in \cite{AGZ}. It
reads (with $k=k_p$ for a $p$ between 1 and $m$)
$$
2k^2 C^2_{\frac{k-1}2}+k^2 C^2_{\frac k2} +\sum_{r=3}^\infty \frac{2k^2}r\left(
\sum_{\Sb k_i\ge 0\\ 2\sum_{i=1}^r k_i=k-r\endSb}\prod_{i=1}^r C_{k_i}\right)^2,
\tag 4
$$
where $\{C_k\}_{k\ge 1}$ are the Catalan numbers, and we assume $C_a=0$ unless
$a\in\{0,1,2,\dots\}$. The Catalan number $C_k$ counts the number of rooted planar
trees with $k$ edges, and different terms of \thetag{4} have the following
interpretation (see \cite{AGZ} for detailed explanations):

$\bullet$ The first term comes from two trees with $(k-1)/2$ edges each that hang
from a common vertex; the factor $k^2$ originates from choices of certain
starting points on each tree united with the common vertex, and the extra 2 is
actually $\E Y_1^2$.

$\bullet$ The second term comes from two trees with $k/2$ edges each
that are glued along one edge. There are $k/2$ choices of this edge for each
of the trees, there is an additional $2=\E Z_{12}^4-1$, and another addional
2 responsible of the choice of the orientation of the gluing.

$\bullet$ The third term comes from two graphs each of which is a cycle of length
$r$ with pendant trees hanging off each of the vertices of the cycle; the total number
of edges in the extra trees being $(k-r)/2$ (this must be an integer). As for the
first term, there is an extra $k^2=k\cdot k$ coming from the choice of the starting
points and also an extra 2 for the choice of the gluing orientation along the cycle.

For each of the three terms the total number of vertices in the resulting graph is equal 
to $k$,
and if one labels each vertex with a letter from an alphabet of cardinality $|B|$
this would yield a factor of
$$
|B|(|B|-1)\cdots (|B|-k+1)=|B|^k+O(|B|^{k-1}).
$$
Normalization by $|B|^k$ yields \thetag{4}.

In the general case, in order to evaluate the covariance
$$
L^{-\frac {k_p+k_q}2}\E\left[\left(\tr\bigl(X_{B_p}^{k_p}(t_p)\bigr)-
\E \tr\bigl(X_{B_p}^{k_p}(t_p)\bigr)\right)\left(\tr\bigl(X_{B_q}^{k_q}(t_q)\bigr)-
\E \tr\bigl(X_{B_q}^{k_q}(t_q)\bigr)\right)\right]
\tag 5
$$
in the limit, we need to employ the same graph counting, except for the two graphs
being glued now correspond to different values $k_p$ and $k_q$ of $k$, and their
vertices are marked by letters of different alphabets $B_p$ and $B_q$. 

$\bullet$ The first term gives
$2k_pk_q C_{\frac{k_p-1}2}C_{\frac{k_q-1}2}$ for the graph counting, and an extra
$$
|B_p\cap B_q|\cdot(|B_p|-1)(|B_p|-2)\cdots (|B_p|-\tfrac{k_p+1}2)\cdot(|B_q|-1)(|B_q|-2)\cdots
(|B_q|-\tfrac{k_q+1}2)
$$
for the vertex labeling (the factor $|B_p\cap B_q|$ comes from the only
common vertex). Moreover, $\E Y_1^2$ is replaced by $\E Y_1(t_p)Y_1(t_q)=c(t_p,t_q)$. 
Normalized by $L^{-\frac {k_p+k_q}2}$ this yields
$$
2k_pk_q C_{\frac{k_p-1}2}C_{\frac{k_q-1}2} c(t_p,t_q)b_{pq}b_p^{\frac{k_p-1}{2}}b_q^{\frac{k_q-1}{2}}.
$$

$\bullet$ The second term has $k_pk_qC_{\frac {k_p}2} C_{\frac {k_q}2}$ from the graph counting
and $c_{pq}^2b_p^{\frac{k_p}{2}-1}b_q^{\frac{k_q}{2}-1}$ from the label counting. In
addition, $\E Z_12^4-1$ is replaced by $\E Z_{12}^2(t_p)Z_{12}^2(t_q)-1=2c^2(t_p,t_q)$.
The total contribution is
$$
k_pk_qC_{\frac {k_p}2} C_{\frac {k_q}2}\bigl(c(t_p,t_q)b_{pq}\bigr)^2
b_p^{\frac{k_p}{2}-1}b_q^{\frac{k_q}{2}-1}.
$$

$\bullet$ For the third term in the same way we obtain
$$
\sum_{r=3}^\infty \frac{2k_pk_q}r\left(
\sum_{\Sb s_i\ge 0\\ 2\sum_{i=1}^r s_i=k_p-r\endSb}\prod_{i=1}^r C_{s_i}\right)
\left(\sum_{\Sb t_i\ge 0\\ 2\sum_{i=1}^r t_i=k_q-r\endSb}\prod_{i=1}^r C_{t_i}\right)
\bigl(c(t_p,t_q)b_{pq}\bigr)^r b_p^{\frac{k_p-r}2} b_q^{\frac{k_q-r}2}
$$
where $c^r(t_p,t_q)$ appeared as $(\E Z_{12}(t_p) Z_{12}(t_q))^r$, which in its turn
came from the edges of the $r$-cycle.

Thus, the asymptotic value of the covariance \thetag{5} is
$$
\multline
2k_pk_q C_{\frac{k_p-1}2}C_{\frac{k_q-1}2} \bigl(c(t_p,t_q)b_{pq}\bigr)b_p^{\frac{k_p-1}{2}}b_q^{\frac{k_q-1}{2}}+
k_pk_qC_{\frac {k_p}2} C_{\frac {k_q}2}\bigl(c(t_p,t_q)b_{pq}\bigr)^2b_p^{\frac{k_p}{2}-1}b_q^{\frac{k_q}{2}-1}
\\+\sum_{r=3}^\infty \frac{2k_pk_q}r\left(
\sum_{\Sb s_i\ge 0\\ 2\sum_{i=1}^r s_i=k_p-r\endSb}\prod_{i=1}^r C_{s_i}\right)
\left(\sum_{\Sb t_i\ge 0\\ 2\sum_{i=1}^r t_i=k_q-r\endSb}\prod_{i=1}^r C_{t_i}\right)
\bigl(c(t_p,t_q)b_{pq}\bigr)^r b_p^{\frac{k_p-r}2} b_q^{\frac{k_q-r}2}.
\endmultline
$$

We now use the fact that for any $S=0,1,2,\dots$
$$
\sum_{\Sb s_i\ge 0\\ \sum_{i=1}^r s_i=S\endSb}\prod_{i=1}^r
C_{s_i}=\binom{2S+r}{S}\frac{r}{2S+r},
$$
see (5.70) in \cite{GKP}. This allows us to rewrite the asymptotic
covariance in terms of binomial coefficients:
$$
\multline
2\binom{k_p}{(k_p-1)/2}\binom{k_q}{(k_q-1)/2}\bigl(c(t_p,t_q)b_{pq}\bigr)b_p^{\frac{k_p-1}{2}}b_q^{\frac{k_q-1}{2}}
\\+
4\binom{k_p}{k_p/2-1}\binom{k_q}{k_q/2-1}
\bigl(c(t_p,t_q)b_{pq}\bigr)^2b_p^{\frac{k_p-2}{2}}b_q^{\frac{k_q-2}{2}}
\\+\sum_{r=3}^\infty 2r \binom{k_p}{(k_p-r)/2}\binom{k_q}{(k_q-r)/2}
\bigl(c(t_p,t_q)b_{pq}\bigr)^r b_p^{\frac{k_p-r}2} b_q^{\frac{k_q-r}2}
\\=\sum_{r=1}^\infty 2r \binom{k_p}{(k_p-r)/2}\binom{k_q}{(k_q-r)/2}
\bigl(c(t_p,t_q)b_{pq}\bigr)^r b_p^{\frac{k_p-r}2} b_q^{\frac{k_q-r}2}
\endmultline
$$

Using the binomial theorem, we can write this expression as a double contour
integral
$$
\frac 2{(2\pi
i)^2}\iint\limits_{const_1=|z|<|w|=const_2}\left(z+\frac{b_p}z\right)^{k_p}
\left(w+\frac{b_q}w\right)^{k_q}\frac{c(t_p,t_q)b_{pq}}{b_p}\frac{dzdw}
{\bigl(\frac{c(t_p,t_q)b_{pq}}{b_p}z-w\bigr)^2}\,.
\tag 6
$$

Consider the right-hand side of \thetag{3} and assume that $|z|^2=b_p<b_q=|w|^2$.
Observe that
$$
\multline
2\ln\left|\frac{c(t_p,t_q)b_{pq}-zw}
{c(t_p,t_q)b_{pq}-z\overline{w}}\right|=-2\ln\left|\frac{\frac{c(t_p,t_q)b_{pq}}{b_p}z-w}
{\frac{c(t_p,t_q)b_{pq}}{b_p}\overline{z}-\overline{w}}\right|
\\=-\ln\left(\frac{c(t_p,t_q)b_{pq}}{b_p}\,z-w\right)
+\ln\left(\frac{c(t_p,t_q)b_{pq}}{b_p}\,z-\overline{w}\right)\\+\ln\left(\frac{c(t_p,t_q)b_{pq}}{b_p}
\,\overline{z}-{w}\right)
-\ln\left(\frac{c(t_p,t_q)b_{pq}}{b_p}\,\overline{z}-\overline{w}\right).
\endmultline
$$
This allows us to rewrite the right-hand side of \thetag{3} as a double contour integral over
complete circles in the form
$$
-\frac{k_pk_q}{2\beta \pi^2}
\oint\limits_{|z|^2=b_p} \oint\limits_{|w|^2=b_q}
(x(z))^{k_p-1}(x(w))^{k_q-1}
\ln\left(\frac{c(t_p,t_q)b_{pq}}{b_p}z-w\right)
\frac{dx(z)}{dz}\frac{dx(w)}{dw}\,dzdw.
$$

Recalling that $\beta=1$ and noting that
$$
k_p(x(z))^{k_p-1}\frac{dx(z)}{dz}=\frac{d(x(z))^{k_p}}{dz}\,,\quad
k_q(x(w))^{k_q-1}\frac{dx(w)}{dw}=\frac{d(x(w))^{k_q}}{dw}\,,
$$
we integrate by parts in $z$ and $w$ and recover \thetag{6}. The proof for for
$b_p=b_q$ is obtained by continuity of both sides, and to see that the needed
identity holds for $b_p>b_q$ it suffices to observe that both sides are symmetric in
$p$ and $q$.

The argument in the case of Hermitian Wigner matrices is exactly the same, except
in the combinatorial part for the first term the factor 2 is missing due to the
change in $\E Y_1(s)Y_1(t)$, in the second term 2 is missing due to the change in $\E
|Z_{12}(s)|^2|Z_{12}(t)|^2$, and in the third term 2 is missing because there is no 
choice in the
orientation of two $r$-cycles that are being glued together. \qed

\subhead Proof of Proposition 1
\endsubhead We need to show that for any complex numbers $\{u_k\}_{k=1}^M$
$$
\sum_{k,l=1}^M u_k\overline{u_l}\int_\H\int_\H f_k(z,s)f_l(w,t) C(z,s;w,t)
\,dzd\bar{z}ds\,dwd\bar{w}dt\ge 0.
$$
We can approximate the integration over
the three-dimensional domains by finite sums of one-dimensional
integrals over semi-circles of the form $|z|=\operatorname{const}$,
$s=\operatorname{const}$. On each 
semi-circle we further uniformly approximate the (continuous) integrand by a polynomial
in $\Re(z)$. Finally, for the polynomials the nonnegativity follows from
Theorem 2'.\qed

\subhead Chebyshev polynomials
\endsubhead One way to describe the limiting covariance structure in the one-matrix
static case is to show that traces of the Chebyshev polynomials of the matrix are 
asymptotically independent, see \cite{J}. A similar effect takes place for
time-dependent submatrices as well.

For $n=0,1,2,\dots$ let $T_n(x)$ be the $n$th degree Chebyshev polynomial of the first kind:
$$
T_n(x)=\cos(n\arccos x),\qquad T_n(\cos(x))=\cos(nx).
$$
For any $a>0$, let $T_n^a(x)=T_n(\frac xa)$ be the rescaled version of $T_n$.

\proclaim{Proposition 3} In the assumptions of Theorem 2'', for any $p,q=1,\dots,m$
$$
\multline
\lim_{L\to\infty}\E\Biggl[\left(\tr\bigl(T_{k_p}^{2\sqrt{b_pL^{k_p}}}(X_{B_p}(t_p))\bigr)-
\E \tr\bigl(T_{k_p}^{2\sqrt{b_pL^{k_p}}}(X_{B_p}(t_p))\bigr)\right)\\ 
\times \left(\tr\bigl(T_{k_q}^{2\sqrt{b_qL^{k_q}}}(X_{B_q}(t_q)))\bigr)-
\E \tr\bigl(T_{k_p}^{2\sqrt{b_pL^{k_q}}}(X_{B_q}(t_q))\bigr)\right)\Biggr]
\\=\delta_{k_pk_q}\,\frac{k_p}{2\beta} \left(\frac{c(t_p,t_q)b_{pq}}{\sqrt{b_pb_q}}\right)^{k_p}.
\endmultline
$$
\endproclaim
\demo{Proof}
Using \thetag{6}  and assuming $b_p<b_q$ we obtain that the needed limit equals
$$
\multline
\dfrac 2{\beta(2\pi
i)^2}\displaystyle\iint\limits_{b_p=|z|<|w|=b_q}T_{k_p}(\cos(\arg(z))
T_{k_q}(\cos(\arg(w))\dfrac{c(t_p,t_q)b_{pq}}{b_p}\dfrac{dzdw}
{\bigl(\frac{c(t_p,t_q)b_{pq}}{b_p}z-w\bigr)^2}\\
=\dfrac 1{2\beta(2\pi
i)^2}\displaystyle\iint\limits_{b_p=|z|<|w|=b_q}\Biggl(\biggl(\frac{z}{\sqrt{b_p}}\biggr)^{k_p}+
\biggl(\frac{\sqrt{b_p}}{z}\biggr)^{k_p}\Biggr) \Biggl(\biggl(\frac{w}{\sqrt{b_q}}\biggr)^{k_q}+
\biggl(\frac{\sqrt{b_q}}{w}\biggr)^{k_q}\Biggr)\\ \times
\dfrac{c(t_p,t_q)b_{pq}}{b_p}\dfrac{dzdw}{\bigl(\frac{c(t_p,t_q)b_{pq}}{b_p}z-w\bigr)^2}.
\endmultline
$$
Writing $(\frac{c_{pq}}{b_p}z-w)^{-2}$ as a series in $z/w$ we arrive at the result. 
Continuity and symmetry of both sides of the limiting relation removes the assumption $b_p<b_q$.\qed
\enddemo

Note that in the Gaussian specialization (when $c(s,t)=\exp(-|s-t|)$) and for 
a single size $L$ time-dependent Wigner matrix (i.e. $b_p=b_q=b_{pq}=1$), 
the centralized traces of Chebyshev polynomials of
this matrix evolve as independent Ornstein-Uhlenbeck processes with speeds equal to
the degrees of the polynomials.


\Refs
\widestnumber\key{ANvM2}

\ref\key ANvM1
\by M.~Adler, E.~Nordenstam, P.~van Moerbeke
\paper The Dyson Brownian minor process
\paperinfo Preprint, 2010, {\tt arXiv:1006.2956}
\endref

\ref\key ANvM2
\by M.~Adler, E.~Nordenstam, P.~van Moerbeke
\paper Consecutive Minors for Dyson's Brownian Motions 
\paperinfo Preprint, 2010, {\tt arXiv:1007.0220}
\endref

\ref \key AGZ
\by G.~W.~Anderson, A.~Guionnet, and O.~Zeitouni
\book An introduction to random matrices
\publ Cambridge University Press \yr 2010
\endref

\ref\key B
\by A.~Borodin
\paper CLT for spectra of submatrices of Wigner random matrices
\paperinfo Preprint, 2010, {\tt arXiv:1010.0898}
\endref

\ref\key BF
\by A.~Borodin and P.~L.~Ferrari
\paper Anisotropic growth of random surfaces in 2+1 dimensions
\paperinfo Preprint, 2008, {\tt arXiv:0804.3035}
\endref

\ref\key FF
\by P.~L.~Ferrari and R.~Frings 
\paper On the partial connection between random matrices and 
interacting particle systems
\paperinfo Preprint, 2010, {\tt arXiv:1006.3946}
\endref


\ref\key GKP
\by R.~L.~Graham, D.~E.~Knuth, and O.~Patashnik
\book Concrete mathematics. A foundation for computer science
\publ Addison-Wesley Publishing Company\publaddr Reading, MA\yr 1994
\endref

\ref\key J
\by K.~Johansson
\paper
On fluctuations of eigenvalues of random Hermitian matrices.
\jour Duke Math. J. \vol 91 \yr 1998\issue 1\pages 151--204
\endref

\ref\key S
\by H.~Spohn
\paper Dyson's model of interacting Brownian motions at arbitrary coupling strength
\jour Markov Processes and Related Fields \vol 4\pages 649--662 \yr 1998
\endref

\endRefs
\enddocument

\end